\documentclass[12pt]{article}

\usepackage{amsmath,amssymb}

\newcommand{\Li}{\text{Li}}
\newcommand{\Ll}{\text{L}}
\newcommand{\Rr}{\mathbb{R}}
\newcommand{\DD}{\mathcal{D}}

\newcommand{\RE}{\mathrm{Re}}
\newcommand{\IM}{\mathrm{Im}}

\begin{document}

\title{Volkov's Pentagon for the Modular Quantum Dilogarithm\footnote{Supported
    by the program ``Mathematical problems of nonlinear dynamics''
    of the Russian Academy of Sciences, RFBR-CNRS grant 09-01-93108-NTsNIL\_a,
    and RFBR grant 11-01-00570}}
\author{L.~D.~Faddeev\\
    St.Petersburg department of Steklov Mathematical Institute
}
\date{January 12, 2012}
\maketitle

\thispagestyle{empty}

\begin{abstract}
    The new form of pentagon equations suggested by Volkov for the
$ q $-exponential on the basis of formal series is derived within the
    Hilbert space framework for the modular version of the quantum
    dilogarithm.
\end{abstract}


\tableofcontents

\section{Motivation}
    The function
$ e(x) $
\begin{align}
\label{f1}
    e(x) = & 1+ \sum_{n=1}^{\infty} 
	\frac{q^{\frac{n(n-1)}{2}}}{(q^{-1}-q)\ldots(q^{-n}-q^{n})} x^{n} = \\
\label{f2}
	= & \prod_{n=0}^{\infty} (1+q^{2n+1}x) 
\end{align}
    has been known since Euler.
    The less familiar representation is given by
\begin{equation*}
    e(x) = \exp \sum_{n=1}^{\infty} \frac{x^{n}(-1)^{n}}{n(q^{n}-q^{-n})} .
\end{equation*}
    Here the denominator of each term in the sum contains an integer
$ n $
    and its
$ q $-counterpart, so that it is natural to speak of the
$ q $-dilogarithm, or, following modern fashion, of the quantum dilogarithm.
    This is especially natural if one uses
$ e(x) $
    in a noncommutative framework.

    Let
$ U $
    and
$ V $
    form a Weyl pair,
\begin{equation}
\label{f3}
    UV = q^{2} VU .
\end{equation}
    Sch\"utzenberger
\cite{Sch} showed that
\begin{equation}
\label{f5}
    e(U) e(V) = e(U+V) .
\end{equation}
    Volkov and I
\cite{FV3}
    derived the formula
\begin{equation*}
    e(V) e(U) = e(U+V+q^{-1}UV) ,
\end{equation*}
    which, by virtue of
(\ref{f3}),
    implies the pentagon relation
\begin{equation}
\label{f7}
    e(V) e(U) = e(U) e(q^{-1}UV) e(V) .
\end{equation}
    Kashaev and I
\cite{FK4}
    showed that this formula is a noncommutative analog of five-term relation
\begin{equation}
\label{f8}
    \Ll(x) + \Ll(y) = \Ll\bigl(\frac{x(1-y)}{1-xy}\bigr) + \Ll(xy)
	+ \Ll\bigl(\frac{y(1-x)}{1-xy}\bigr) ,
\end{equation}
    for the Rogers dilogarithm
$ \Ll(x) $,
    where
$ x $ and $ y $
    are variables in the interval
$ [0,1] $.
    The proof is based on the quasiclassical computation in the limit as
$ q \to 1 $
    with the use of the asymptotic
\begin{equation*}
    e(u) = \exp \frac{1}{2\ln q} E(u) , \quad
	E(u) = \Li_{2}(-u) ,
\end{equation*}
    where
$ \Li_{2}(u) $ 
    is the Euler dilogarithm
\begin{equation*}
    \Li_{2}(u) = \sum \frac{u^{n}}{n^{2}}.
\end{equation*}
    The functions
$ \Ll(x) $ and
$ \Li_{2}(x) $
    are related by the formula
\begin{equation*}
    \Ll(x) = \Li_{2}(x) +\frac{1}{2} \ln x \ln(1-x) , \quad 0 \leq x \leq 1.
\end{equation*}

    The last result gives rise to some questions:

    1. Why does the noncommutative formula involve the Euler dilogarithm,
    while the Rogers dilogarithm occurs in the classical limit?

    2. Why are the arguments in the functions in formulas
(\ref{f7}) and
(\ref{f8})
    so drastically different?

    The answer to the first question was already discussed in
\cite{FK4}
    and the definitive answer was found in 
\cite{KN5},
    where the quasiclassical asymptotic in formula
(\ref{f7})
    and its generalizations arising in conformal field theory were studied
    in detail.
    A very beautiful answer to the second question was given by Volkov
\cite{V1},
    who rewrote the pentagon relation in a new form.

    To describe it I shall start from the classical case. Let
\begin{equation*}
    R(u) = \Ll(\frac{u}{1+u}) .
\end{equation*}
    One can readily verify that relation
(\ref{f8})
    is equivalent to
\begin{equation}
\label{f12}
    R(u) + R(v) = R(\frac{v}{1+u}) + R(\frac{uv}{1+u+v}) + R(\frac{u}{1+v}) .
\end{equation}
    The arguments
\begin{equation*}
    x_{1} = u, \quad x_{2} = v, \quad x_{3} = \frac{1+v}{u} , \quad
	x_{4} = \frac{1+u+v}{uv} , \quad x_{5} = \frac{1+u}{v} ,
\end{equation*}
    occurring here, give a solution of the
$ Y $-system
$ A_{2} $,
\begin{equation}
\label{f14}
    x_{i} x_{i+2} = 1+x_{i+1} ,
\end{equation}
    which has period 5,
\begin{equation*}
    x_{i+5} = x_{i} .
\end{equation*}
    Thus
(\ref{f12})
    can be rewritten in the form
\begin{equation}
\label{f16}
    R(x_{1}) + R(x_{2}) = R(\frac{1}{x_{5}}) + R(\frac{1}{x_{4}})
	+ R(\frac{1}{x_{3}}) .
\end{equation}

    The quantum counterpart of system
(\ref{f14})
    has the form
\begin{equation}
\label{f17}
    X_{i} X_{i+2} = 1+q X_{i+1}
\end{equation}
    and if
\begin{equation*}
    X_{1} = U , \quad X_{2} = V , \quad UV=q^{2} VU ,
\end{equation*}
    then
\begin{align}
\nonumber
    X_{3} = & U^{-1} (1+qV) , \\
\label{f19}
    X_{4} = & U^{-1} (q^{-1} + U + V) V^{-1}  \\
\nonumber
    X_{5} = & (1+qU) V^{-1} .
\end{align}
    Volkov's quantum formula reads
\begin{equation}
\label{f20}
    e(X_{1}) e(X_{2}) = e(X_{5}^{-1}) e(X_{4}^{-1}) e(X_{3}^{-1}) .
\end{equation}
    We see a remarkable analogy between formulas
(\ref{f16}) and
(\ref{f20}).

    Volkov's derivation is based on calculations with formal series.
    The nature of the parameter
$ q $ was not material.
    The attempt to embed his argument in a Hilbert space framework
    encounters an obstacle.
    In particular, the simplest involution
\begin{equation*}
    U^{*} = U , \quad V^{*} = V
\end{equation*}
    requires that
$ |q|=1 $,
    but the original formulas
(\ref{f1}) and
(\ref{f2})
    are undefined if
$ \frac{\ln q}{i\pi} $
    is rational.

    To overcome this difficulty I suggested
\cite{F6}
    to replace the function
$ e(u) $
    by its modular counterpart
\begin{equation*}
    \Theta(u) = \frac{\prod(1+q^{2n+1}u)}{\prod(1+\tilde{q}^{2n+1}\tilde{u})} ,
\end{equation*}
    where
\begin{equation*}
    q = e^{i\pi\tau} , \quad \tilde{q} = e^{-i\pi/\tau} ,
	\tilde{u} = u^{1/\tau} .
\end{equation*}

    In the present paper I show that formula
(\ref{f20})
    remains valid if
$ e(u) $ 
    is replaced by
$ \Theta(u) $
    and give a derivation of this formula based on the solution of the
$ Y $-system
(\ref{f17}).

\section{Modular quantum dilogarithm}
    Let
\begin{equation}
\label{f24}
    \gamma(z) = \exp \bigl\{ -\frac{1}{4}\int_{-\infty}^{\infty}
	\frac{e^{itz}}{\sin \omega t \sin \omega' t} \frac{dt}{t} \bigr\},
\end{equation}
    where the integration contour bypasses the singularity at
$ t=0 $
    from above.
    Here, following Volkov
\cite{V7},
    I use the notation
$ \omega $ and $  \omega' $,
    dated to Weierstrass's times, for the periods of the lattice with
    generators 1 and
$ \tau $,
    where
\begin{equation*}
    \tau = \frac{\omega'}{\omega} , \quad \omega\omega' = - \frac{1}{4} .
\end{equation*}
    The integral in
(\ref{f24})
    is well defined for
$ \omega $ and
$ \omega' $
    lying in the upper half-plane, and
$ \omega $ and
$ \omega' $
    corresponding to a positive    
$ \tau $
    are pure imaginary.
    One can readily verify that
\begin{equation}
\label{f26}
    \Theta(e^{-i\pi z/\omega}) = \gamma(z) .
\end{equation}

    We take formula
(\ref{f26})
    to be the definition of the modular quantum dilogarithm.
    The function
$ \gamma(z) $, as well as
$ e(u) $,
    is well known in analysis.
    It has many names; in particular, it  is called the ``double
$ \Gamma $-function'' in
\cite{B8}.
    Its properties, including noncommutative versions, are described
    in detail in Volkov's paper
\cite{V7}
    and in some Kashaev's papers 
\cite{K9}, \cite{K10}.
    I introduced it in
\cite{F6}
    to describe the duality in the discrete Weyl-Heisenberg group and used it
    in
\cite{F7}
    to define the modular double of the quantum group
$ SL_{q}(2) $,
    where it replaces the
$ q $-exponential in the definition of Drinfeld's universal 
$ R $-matrix
\cite{D11}.
    It is used in the papers
\cite{FV12} and \cite{FKV13}
    to describe the evolution operator for quantum integrable systems
    on discrete space-time, and the paper
\cite{FV14}
    gives the definition of fundamental
$ R $-matrix for the XXZ model of arbitrary spin.
    Further applications of this function, in particular, in mathematics
    are listed in the introduction to the paper
\cite{KN5}.

    Let us present the properties of
$ \gamma(z) $,
    to be used in what follows

    1. The functional equation
\begin{equation}
\label{f27}
    \frac{\gamma(z+\omega')}{\gamma(z-\omega')} = 1 +e^{-i\pi z/\omega} .
\end{equation}
    Note that a similar property with
$ \omega $ and
$ \omega' $
    interchanged hold as well, but we do not need it.

    2. The inversion formula
\begin{equation}
\label{f28}
    \gamma(z) \gamma(-z) = e^{i\beta} e^{i\pi z^{2}} , \quad
	\beta = \frac{\pi}{12} (\tau + \frac{1}{\tau}) .
\end{equation}

    3. Unitarity
\begin{equation*}
    \overline{\gamma(z)} = \frac{1}{\gamma(z)}
\end{equation*}
    for
$ \tau >0 $
    and real
$ z $.

    4. The quasiclassical asymptotics
\begin{equation}
\label{f30}
    \gamma(z) = \exp \frac{1}{2\pi i\tau} E(e^{-i\pi z/\omega}) .
\end{equation}

    5. Asymptotics
\begin{equation}
\label{asymp}
    \gamma(z) \to 1, \quad \RE z \to \infty 
\end{equation}
    with bounded
$ \IM z $.
    For
$ \RE z \to -\infty $
    one can use the inversion formula.

    6. Zeros and poles\\
    The zero of
$ \gamma(z) $
    nearest to real axis is at
$ z=\omega'' $, where
\begin{equation*}
    \omega'' = \omega + \omega' .
\end{equation*}
    The pole of
$ \gamma(z) $
    nearest to real axis is at
$    z = - \omega'' $. 
    The coefficients in the main terms
\begin{equation}
\label{mt}
    \gamma(z) \simeq c_{1} (z-\omega'') , \quad
	\gamma(z) \simeq \frac{c_{2}}{z+\omega''}
\end{equation}
    are given by 
\begin{equation*}
    c_{1} = 2\pi ic , \quad c_{2} = -\frac{c}{2\pi i} , \quad
	c = e^{-\frac{i\pi}{12}(\tau+\frac{1}{\tau}) - \frac{i\pi}{4}} .
\end{equation*}

    7. The Fourier transform duality
\begin{equation}
\label{FTD}
    \int_{-\infty}^{\infty} \gamma(t-\omega''+i0) e^{-2\pi i xt} dt
	= c\frac{1}{\gamma(x+\omega''-i0)}
\end{equation}
    for real
$ x $.

    8. The main integral relation
\begin{equation}
\label{MIR}
    \int_{-\infty}^{\infty} \frac{\gamma(t-\omega''+i0)}{\gamma(x+t)}
	e^{-2\pi yt} dt = c\frac{\gamma(x+y)}{\gamma(x)\gamma(y+\omega''-i0)} 
\end{equation}
    for real
$ y $.
    We see, that the pole structure
(\ref{mt})
    and asymptotic behaviour
(\ref{asymp})
    in the LHS and RHS are compatible.
    The value of
$ c $
    can be obtained via the inversion formula
\begin{multline*}
    c^{2} = -c_{1}c_{2} = \lim_{\epsilon \to 0} \gamma(\omega''-i\epsilon)
	\gamma(-\omega''+i\epsilon) = e^{i\beta} e^{i\pi\omega''^{2}} =\\
	= \exp i\pi(-\frac{1}{6}(\tau+\frac{1}{\tau})-\frac{1}{2}) ,
\end{multline*}
    where we used that
\begin{equation*}
    \omega''^{2} = \omega^{2} + \omega'^{2} +2\omega\omega'
	= \omega\omega'(\tau+\frac{1}{\tau}+2)
	= -\frac{1}{4}(\tau+\frac{1}{\tau}+2) .
\end{equation*}
    We can shift the integration contour in
(\ref{FTD}) and
(\ref{MIR}),
    for instance
\begin{equation}
\label{ShC}
    \int_{-\infty}^{\infty}\gamma(t) e^{-2\pi ixt} dt 
	= c e^{2\pi i x\omega''} \frac{1}{\gamma(x+\omega''-i0)}.
\end{equation}

    Formula
(\ref{f27})
    implies the relation
\begin{equation}
\label{f31}
    \frac{\Theta(qu)}{\Theta(q^{-1}u)} = \frac{1}{1+u} .
\end{equation}

\section{The evolution operator}
    Now we are in a position to proceed to the construction of the evolution
    operator for the quantum
$ A_{2} $
    system.
    Let us realize the operators
$ U $ and
$ V $
    on a dense domain in
$ L_{2}(\Rr) $
    by the formulas
\begin{equation*}
    U f(z) = e^{-i\pi z/\omega} f(z) , \quad
    V f(z) = f(z+ 2\omega') .
\end{equation*}
    These formulas involve multiplication by an unbounded function and a shift
    into the complex plane.
    The domain
$ \DD $
    of
$ U $ and
$ V $
    should permit these operations.
    It suffices to assume that
$ \DD $
    contains entire functions of
$ z $,
    decaying at the rate
$ e^{-z^{2}} $
    in the directions parallel to the real axis.
    For
$ \tau >0 $
    the operators
$ U $ and
$ V $
    are nonnegative and essentially selfadjoint.
    The same is true for the operators
$ X_{3} $,
$ X_{4} $ and
$ X_{5} $
    in
(\ref{f19}).

    Let us introduce the operator
$ K $
    as a multiplication by
$ \Theta(e^{-i\pi z/\omega}) $
\begin{equation*}
    K f(z) = \Theta(e^{-i\pi z/\omega}) f(z) = \gamma(z) f(z)
\end{equation*}
    and let 
$ F $
    be the Fourier transform
\begin{equation*}
    F f(z) = \int_{-\infty}^{\infty} e^{-2\pi i zt} f(t) dt .
\end{equation*}
    Note that
\begin{align*}
    UF =& FV , \quad VF = FU^{-1} \\
    V\Theta(U) =& \Theta(U) (1+q^{-1}U) V .
\end{align*}
    The last property follows from the functional equation
(\ref{f27}).
    Now one can readily see that
\begin{equation}
\label{f38}
    X_{i+1} = S^{-1} X_{i} S ,
\end{equation}
    where the unitarity operator
$ S $
    is given by the formula
\begin{equation}
\label{f39}
    S = K F .
\end{equation}
    The periodicity
$ i \equiv i+5 $
    will be insured if
\begin{equation}
\label{f40}
    S^{5} = e^{i\alpha} I .
\end{equation}
    Let us show that this is indeed true. Let
$ S^{5}(x,y) $
    be the integral kernel of the operator
$ S^{5} $.
    By definition
\begin{multline*}
    S^{5}(x,y) = \gamma(x) \int \exp \{-2\pi i (xt_{1} +t_{1}t_{2}
	+t_{2}t_{3} +t_{3}t_{4} +t_{4}y) \} \times\\
    \times \gamma(t_{1}) \gamma(t_{2})
	\gamma(t_{3}) \gamma(t_{4}) dt_{1} dt_{2} dt_{3} dt_{4} .
\end{multline*}
    We shall show, that
\begin{equation*}
    S^{5}(x,y) = e^{i\alpha} \delta(x-y)
\end{equation*}
    and give the explicit value of
$ \alpha $.
    To do it we calculate integrals using formulas
(\ref{FTD})--(\ref{ShC}).
    Let us denote
$ t_{2}=t $
    and
$ t_{4}=s $
    and take integrals over
$ t_{1} $
    and
$ t_{3} $
    to get
\begin{multline*}
    S^{5}(x,y) = c^{2} \int \frac{\gamma(x)\gamma(t)\gamma(s)}{
	\gamma(x+t+\omega'')\gamma(s+t+\omega'')} \times \\
    \times \exp\bigl\{2\pi i\bigl((x+t)\omega'' +(s+t)\omega'' -sy \bigr)
	\bigr\} dt \,ds.
\end{multline*}
    Then integration over
$ s $
    gives
\begin{equation*}
    S^{5}(x,y) = c^{3} \frac{\gamma(x)}{\gamma(y)}
	\int \frac{\gamma(t+y-\omega'')}{\gamma(t+x+\omega'')}
    \exp \bigl\{2\pi i(2t\omega'' +x\omega'' +y\omega''-\omega''^{2} )
	\bigr\} dt .
\end{equation*}
    and integral over
$ t $
    gives
\begin{equation*}
    S^{5}(x,y) = c^{4} \frac{\gamma(x)}{\gamma(y)} \exp \bigl\{
	2\pi i\bigl((x-y)\omega'' -\omega''^{2}\bigr)
    \bigr\} \frac{\gamma(x-y-\omega''+i0)}{\gamma(x-y+\omega''-i0)
	\gamma(-\omega''+i0)} .
\end{equation*}
    Here I restored the rule of going around the zeros and poles of
$ \gamma(z) $
    by writing
$ \omega''-i0 $
    instead of
$ \omega'' $.

    The expression
\begin{equation*}
    \lim_{\epsilon\to0}
\frac{\gamma(a-\omega''+i\epsilon)}{\gamma(a+\omega''-i\epsilon)}
    \cdot \frac{1}{\gamma(-\omega''+i\epsilon)} = I(a)
\end{equation*}
    vanishes for
$ a \neq 0 $
    and in the vicinity of
$ a = 0 $
    we can substitute the full
$ \gamma(z) $
    by its main terms
(\ref{mt}),
    to get
\begin{equation*}
    I(a) \simeq \frac{c_{2}}{c_{1}c_{2}}
	\frac{i\epsilon}{(a-i\epsilon)(a+i\epsilon)} 
    = \frac{1}{c} \frac{1}{2\pi} \frac{\epsilon}{a^{2}+\epsilon^{2}}
	\to \frac{1}{c} \delta(a) 
\end{equation*}
    and so that
\begin{equation*}
    S^{5}(x,y) = e^{i\alpha} \delta(x-y) .
\end{equation*}
    Thus
(\ref{f40})
    is derived and
\begin{equation*}
    e^{i\alpha} = c^{3} e^{-2\pi i \omega''^{2}} ,
\end{equation*}
    and
\begin{equation}
\label{f43}
    \alpha = \pi\bigl( -\frac{1}{4}(\tau+\frac{1}{\tau})-\frac{3}{4}
	+\frac{1}{2}(\tau+\frac{1}{\tau} +2)\bigr)
    = \frac{i\pi}{4}(\tau+\frac{1}{\tau}+1) = 3\beta +\frac{\pi}{4} .
\end{equation}
    Now we can proceed to the derivation of Volkov's formula.

\section{The main derivation}
    First note that, by virtue of the inversion formula
(\ref{f28}),
    we can give an alternative form of the evolution operator
\begin{equation}
\label{f45}
    S = e^{i\beta} \hat{K}^{-1} G , \quad \beta
	= \frac{\pi}{12}(\tau+\frac{1}{\tau}) ,
\end{equation}
    where
\begin{equation*}
    \hat{K} f(z) = \Theta(e^{i\pi z/\omega}) f(z) = \gamma(-z) f(z)
\end{equation*}
    and
\begin{equation*}
    G(x,y) = e^{i\pi z^{2}} F .
\end{equation*}
    One can readily verify that
\begin{equation}
\label{f48}
    G^{3} = e^{i\pi/4} F^{2} .
\end{equation}

    Now we have all we need to prove Volkov's formula in the form
\begin{equation}
\label{f49}
    \Theta(X_{1}) \Theta(X_{2}) = \Theta(X_{5}^{-1}) \Theta(X_{4}^{-1})
	\Theta(X_{3}^{-1}) .
\end{equation}
    We use the action
(\ref{f38})
    of the operator
$ S $
    to rewrite this relation in the form
\begin{equation}
\label{MR}
    K S^{-1} K S = S^{-4} \hat{K} S^{4} S^{-3} \hat{K}
	S^{3} S^{-2} \hat{K} S^{2} 
\end{equation}
    or in view of definitions
(\ref{f39}) and
(\ref{f45})
\begin{equation*}
    K F^{-1} = S^{-4} G^{3} e^{3i\beta} = e^{-i\alpha +3i\beta}
	K F G^{3} .
\end{equation*}
    Thus relation
(\ref{f49})
    is true by virtue of the definition of
$ \alpha $ and
$ \beta $
    in formulas
(\ref{f43}) and
(\ref{f28}),
    formula
(\ref{f48})
    and the identity
\begin{equation*}
    F^{4} = I .
\end{equation*}
    Thus we have proved the main result
(\ref{f49})
    and shown that it is equivalent to
(\ref{f40}).

\section{Conclusion}
    In conclusion let us verify that formula
(\ref{f40})
    gives the classical identity for the Rogers dilogarithm.
    We proceed to the quasiclassical computation of the kernel
$ S^{5}(x,y) $.
    Let us rewrite it in the form
\begin{equation*}
    S^{5}(x,y) = \int M(x,z) e^{2\pi i(x-y)z} dz
\end{equation*}
    with
\begin{multline*}
    M(x,z) = \gamma(x)\gamma(z) \int\exp \{-2\pi(xt_{1}+t_{1}t_{2}
	+t_{2}t_{3}+t_{3}z+xz)\} \times \\
    \times \gamma(t_{1})\gamma(t_{2})\gamma(t_{3})  dt_{1} dt_{2} dt_{3} .
\end{multline*}
    We have already proved, that
$ M(x,z) $
    does not depend on
$ x $
    and
$ z $
\begin{equation*}
    M(x,z) = e^{i\alpha} .
\end{equation*}
    Let us check it in the quasiclassical approximation.
    By using the asymptotic
(\ref{f30})
    and by retaining only the coefficients of
$ \frac{1}{2\pi i\tau} $
    in the exponential in the integrand, we arrive at the integral
\begin{equation*}
    M(x,y) \sim \int \exp \frac{1}{2\pi i\tau} \bigl\{ \sum_{i=1}^{5}
	(E(e^{p_{i}}) +p_{i}p_{i+1})\bigr\} dp_{3} dp_{4} dp_{5},
\end{equation*}
    where
\begin{equation*}
    p_{1} = - \frac{i\pi x}{\omega}
\end{equation*}
    and the arguments
$ p_{2}, \ldots p_{5} $
    arise in a similar way from
$ y, t_{1}, t_{2} $ and $ t_{3} $, while
$ p_{6} = p_{1} $.
    By using the relation
\begin{equation*}
    \frac{d}{dp} E(e^{p}) = -\ln(1+e^{p})
\end{equation*}
    we obtain the stationary phase equations
\begin{align*}
    \ln(1+e^{p_{3}}) =& p_{2} + p_{4} , \\
    \ln(1+e^{p_{4}}) =& p_{3} + p_{5} , \\
    \ln(1+e^{p_{5}}) =& p_{4} + p_{1} , 
\end{align*}
    which coincide with the three equations in the
$ Y $-system for
$ x_{i} = e^{p_{i}} $,
    expressing
$ x_{3} $, 
$ x_{4} $ and
$ x_{5} $
    via
$ x_{1} = x $ and
$ x_{2} = y $.

    Now by rewriting the sum
$ \sum_{i} p_{i} p_{i+1} $
    in the form
$ \sum_{i} \frac{1}{2} p_{i}(p_{i-1} + p_{i+1}) $
    we obtain the quasiclassical result
\begin{equation*}
    \sum \bigl(E(x_{i}) +\frac{1}{2} \ln x_{i} \ln(1+x_{i}) \bigr)
	= -\frac{\pi^{2}}{2} = -3\frac{\pi^{2}}{6} ,
\end{equation*}
    which leads to formula
(\ref{f16})
    in view of the fact that
\begin{equation*}
    R(x) = -E(x) -\frac{1}{2} \ln x \ln (1+x) 
\end{equation*}
    and
\begin{equation*}
    \frac{\pi^{2}}{6} - R(x) = R(\frac{1}{x}) .
\end{equation*}
    We see explicitly how the Euler dilogarithm gets replaced by the Rogers
    dilogarithm.
    Using the terminology of the classical mechanics and quasiclassical
    experience we can say, that the Euler and Rogers dilogarithms
    give truncated action and full action, correspondingly.

    A majority of results in this paper arose in the course of discussion
    with A.~Yu.~Volkov and R.~M.~Kashaev.
    Some of the results presented here can be found in the literature
    dealing with applications of the quantum dilogarithm; e.~g. see
\cite{ChF15} and \cite{G16}.
    I have just gathered these results into one of possible logical
    schemes. I am grateful to Volkov and Kashaev for participating
    in the discussions. The paper was finished in Greece in the excellent
    conditions provided for me by Professors Cotsiolis, Savvidy
    and Zoupanos. I am keenly grateful to all of them.

\end{document}